\documentclass[11pt]{amsart}
\usepackage[T1]{fontenc}
\usepackage{amsmath,amssymb,amsthm}
\usepackage{microtype}
\usepackage[hidelinks]{hyperref}
\newtheorem{theorem}{Theorem}
\newtheorem{corollary}[theorem]{Corollary}
\newtheorem{proposition}[theorem]{Proposition}
\newtheorem{lemma}[theorem]{Lemma}
\newcommand{\Sch}{\mathsf{S}}
\newcommand{\SubSch}{\mathsf{S}_{\mathrm{sub}}}
\newcommand{\StrSch}{\mathsf{S}_{\mathrm{str}}}
\newcommand{\ca}{\operatorname{ca}}
\newcommand{\tca}{\widetilde{\operatorname{ca}}}
\newcommand{\tdelta}{\widetilde{\delta}}
\newcommand{\K}{\mathbb K}
\newcommand{\N}{\mathbb N}
\title[Quantitative Schur properties]{Sharp constants in the quantitative Schur property}
\author[D. N\'u\~nez-Alarc\'on]{Daniel N\'u\~nez-Alarc\'on}
\address{Departamento de Matem\'aticas, Universidad Nacional de Colombia, Edificio 404, Carrera 30 No. 45-03, Bogot\'a D.C., Colombia}
\email{dnuneza@unal.edu.co}
\author[D. Pellegrino]{Daniel Pellegrino}
\address{Departamento de Matem\'atica, Universidade Federal da Para\'iba, 58051-900 Jo\~ao Pessoa, PB, Brazil}
\email{daniel.pellegrino@academico.ufpb.br}
\author[D. M. Serrano-Rodr\'iguez]{Diana Marcela Serrano-Rodr\'iguez}
\address{Departamento de Matem\'aticas, Universidad Nacional de Colombia, Edificio 404, Carrera 30 No. 45-03, Bogot\'a D.C., Colombia}
\email{diserranor@unal.edu.co}
\begin{document}
\begin{abstract}
We answer a question of C\'uth and Kalenda by showing that, for every $c\ge1$, every $c$-strong Schur Banach space is $(2c+1)$-Schur and that the constant $2c+1$ is optimal. Over either scalar field, we determine the full range of pairs of optimal constants for the quantitative Schur inequality and its subsequence version. The admissible pairs $(a,b)$ are exactly those satisfying $1\le a\le b\le2a+1$, and each is realized by an equivalent norm on $\ell_1$. We also answer a question of Kalenda concerning a sufficient condition for the $1$-Schur property.
\end{abstract}
\maketitle
\section{Introduction}
In 1921, Schur \cite{Schur} proved that weak convergence of sequences in the space of absolutely summable sequences implies norm convergence. A Banach space with this property is said to have the \emph{Schur property}. Throughout, $\K=\mathbb R$ or $\mathbb C$, $\N=\{1,2,\ldots\}$, and $\ell_1(\K)$ denotes the space of absolutely summable scalar sequences, with norm $\|x\|_1=\sum_{i=1}^\infty|x_i|$.

In their study of subspaces of $L_1$ which embed into $\ell_1$, Godefroy, Kalton and Li \cite{GKL} introduced the $1$-strong Schur property. This property gives a quantitative lower $\ell_1$-estimate on subsequences of separated sequences. For $c\ge1$, a Banach space $X$ over $\K$ has the \emph{$c$-strong Schur property} if, for every $d\in(0,2]$ and every $\varepsilon>0$, each sequence $(z_n)_{n\in\N}$ in $X$ satisfying
\[
 \|z_n\|_X=1\quad(n\in\N),\qquad
 \|z_n-z_m\|_X\ge d\quad(n\ne m)
\]
has a subsequence $(z_{n_k})_{k\in\N}$ such that
\begin{equation}\label{eq:strong}
 \frac{1}{2c/d+\varepsilon}\sum_{k=1}^N|\lambda_k|
 \le\left\|\sum_{k=1}^N\lambda_k z_{n_k}\right\|_X
 \le\sum_{k=1}^N|\lambda_k|
\end{equation}
for every $N\in\N$ and $(\lambda_k)_{k=1}^N\in\K^N$ (see \cite[Proposition 5.6]{K}).

For a Banach space $X$ over $\K$, let $X^*$ denote the continuous dual of $X$, $B_{X^*}$ the closed unit ball of $X^*$, and $\ell_\infty(X)$ the Banach space of bounded sequences in $X$, endowed with the supremum norm. For $z=(z_n)_{n\in\N}\in\ell_\infty(X)$, define
\begin{align}
\ca_X(z)&=\inf_{N\in\N}\sup_{j,k\ge N}\|z_j-z_k\|_X,\label{eq:definitionca}\\
\delta_X(z)&=\sup_{f\in B_{X^*}}\ca_{\K}((f(z_n))_{n\in\N}),\label{eq:definitiondelta}
\end{align}
where $\K$ is equipped with its usual absolute value. Thus $\ca_X(z)$ and $\delta_X(z)$ measure, respectively, the failure of $z$ to be norm Cauchy and to be weakly Cauchy. For $c\ge1$, $X$ has the \emph{$c$-Schur property} if
\[
 \ca_X(z)\le c\delta_X(z)\qquad(z\in\ell_\infty(X)).
\]
Kalenda and Spurn\'y \cite{KS} proved that $\ell_1(\K)$ has the $1$-Schur property; in this case, $\ca_X(z)=\delta_X(z)$ for every $z\in\ell_\infty(X)$ with $X=\ell_1(\K)$.

The corresponding quantities for subsequences are
\begin{align}
\tca_X(z)&=\inf_{(n_k)}\ca_X((z_{n_k})_{k\in\N}),\label{eq:definitiontca}\\
\tdelta_X(z)&=\inf_{(n_k)}\delta_X((z_{n_k})_{k\in\N}),\label{eq:definitiontdelta}
\end{align}
where both infima are taken over all strictly increasing sequences $(n_k)_{k\in\N}$ in $\N$. We omit the subscript $X$ when the norm is clear from the context. For every $z\in\ell_\infty(X)$,
\begin{equation}\label{eq:basic}
\delta_X(z)\le\ca_X(z),\qquad \tca_X(z)\le\ca_X(z),\qquad \tdelta_X(z)\le\delta_X(z).
\end{equation}

For every $c\ge1$, Kalenda \cite[Proposition 5.1]{K} proved that the inequalities
\[
 \tca_X(z)\le c\delta_X(z)\qquad(z\in\ell_\infty(X))
\]
imply
\[
 \ca_X(z)\le(2c+1)\delta_X(z)\qquad(z\in\ell_\infty(X)).
\]
He asked whether $2c+1$ is optimal \cite[Question 9.2]{K}. His Examples 8.3 and 8.4 give a lower bound of $2$ when $c=1$, over either scalar field. In the real case, the bound $3$ for $c=1$ improves the earlier bound $5$ in the passage from the $1$-strong Schur property to the quantitative Schur property (see \cite[Remark 5.7]{K}).

In their recent work on Lipschitz-free spaces, C\'uth and Kalenda \cite[Proposition 2.2(c)]{CK} recalled that every real Banach space with the $c$-strong Schur property has the $(2c+1)$-Schur property. They asked whether $2c+1$ is optimal or can be reduced to $2c$ or $c+1$ (see the paragraph following Proposition 2.2 in the preprint \cite{CK}).

For an infinite-dimensional Banach space $X$, write
\begin{align*}
 \Sch(X)&=\inf\{c\ge1:\ca_X(z)\le c\delta_X(z)\text{ for all }z\in\ell_\infty(X)\},\\
 \SubSch(X)&=\inf\{c\ge1:\tca_X(z)\le c\delta_X(z)\text{ for all }z\in\ell_\infty(X)\},\\
 \StrSch(X)&=\inf\{c\ge1:X\text{ has the }c\text{-strong Schur property}\},
\end{align*}
with $\inf\varnothing=+\infty$. The first two infima, when finite, are attained. The same holds for the third: for given $d$ and $\varepsilon$ in \eqref{eq:strong}, take an admissible $c<\StrSch(X)+d\varepsilon/4$ and use \eqref{eq:strong} with tolerance $\varepsilon/2$. The constant $\SubSch(X)$ is unchanged if $\delta_X$ is replaced by $\tdelta_X$: applying the defining inequality to each subsequence $w$ of $z$ gives $\tca_X(z)\le\tca_X(w)\le c\delta_X(w)$, and the converse follows from \eqref{eq:basic}. The restriction $c\ge1$ in the first two infima does not affect their values in infinite dimension. Indeed, Riesz's lemma gives a separated sequence $z$ in the unit sphere, so $d_0=\tca_X(z)>0$. For every $\eta>0$ there is a subsequence $w$ with $\ca_X(w)<d_0+\eta$ and $\tca_X(w)\ge d_0$. Since $\delta_X(w)\le\ca_X(w)$, any admissible constant is at least $d_0/(d_0+\eta)$ for every $\eta>0$.

Kalenda's results \cite[Propositions 5.1 and 5.6]{K} give
\begin{equation}\label{eq:kalenda}
 \SubSch(X)\le\Sch(X)\le2\SubSch(X)+1,
\end{equation}
Moreover, $\StrSch(X)=\SubSch(X)$ if $\K=\mathbb R$, while
\[
 \SubSch(X)\le\StrSch(X)\le\frac{\pi}{2}\SubSch(X)
 \qquad(\K=\mathbb C).
\]
We show that every pair allowed by \eqref{eq:kalenda} occurs on an equivalent renorming of $\ell_1(\K)$, with the strong Schur constant equal to the subsequence constant.

For isomorphic Banach spaces $X$ and $Y$ over $\K$, their Banach--Mazur distance is
\[
 d(X,Y)=\inf\{\|T\|\|T^{-1}\|:\ T:X\to Y\text{ is a bounded linear isomorphism}\}.
\]

\begin{theorem}\label{thm:main}
Let $\K=\mathbb R$ or $\mathbb C$. The set of finite pairs $(\SubSch(X),\Sch(X))$ attained by infinite-dimensional Banach spaces over $\K$ is exactly
\[
 \{(a,b)\in\mathbb R^2:1\le a\le b\le 2a+1\}.
\]
For every pair $(a,b)$ in this region there is an equivalent norm on $\ell_1(\K)$ for which the resulting space $X$ satisfies
\[
 \SubSch(X)=\StrSch(X)=a,\qquad \Sch(X)=d(X,\ell_1(\K))=b.
\]
\end{theorem}

In \cite[Proposition 6.4]{K}, Kalenda proved that a Banach space $X$ over $\K$ has the $1$-Schur and $1$-strong Schur properties if every sequence $(z_n)_{n\in\N}\in\ell_\infty(X)$ admits a subsequence $(z_{n_k})_k$ and a vector $x\in X$ such that for every further subsequence $(z_{n_{k_l}})_l$ and every $y\in\operatorname{span}(\{x\}\cup\{z_{n_k}:k\in\N\})$,
\begin{equation}\label{eq:add}
 \limsup_l\|z_{n_{k_l}}+y\|=\|x+y\|+\limsup_l\|z_{n_{k_l}}-x\|.
\end{equation}
Kalenda showed that this condition is not necessary over the complex field \cite[Remark 6.7(3)]{K} and asked whether it is necessary over the real field \cite[Question 9.3]{K}. The answer is negative: for each $1<p<\infty$, the space $\K\oplus_p\ell_1(\K)$ has both the $1$-Schur and $1$-strong Schur properties but fails the condition in \eqref{eq:add} (Proposition~\ref{prop:counter}). Here $\K\oplus_p\ell_1(\K)$ denotes the space of pairs $(\alpha,u)\in\K\times\ell_1(\K)$ endowed with the norm $(|\alpha|^p+\|u\|_1^p)^{1/p}$.

\section{Proof of Theorem \ref{thm:main}}
We use the following elementary block selection in $\ell_1$.

\begin{lemma}\label{lem:block}
Every bounded sequence in $\ell_1(\K)$ has a subsequence $(v_k)$ for which there exist $x\in\ell_1(\K)$, finitely supported vectors $u_k$ with pairwise disjoint supports, $L\ge0$ and $\sigma\in\K$ such that
\begin{equation}\label{eq:blocks}
 \|v_k-x-u_k\|_1\to0,\qquad \|u_k\|_1\to L,\qquad \sum_i(u_k)_i\to\sigma.
\end{equation}
For every $\eta>0$, after deleting finitely many terms,
\begin{equation}\label{eq:blocklower}
 \left\|\sum_{k=1}^N\lambda_k v_k\right\|_1\ge (L-\eta)\sum_{k=1}^N|\lambda_k|
\end{equation}
for all $N\in\N$ and $\lambda_1,\ldots,\lambda_N\in\K$.
\end{lemma}
\begin{proof}
By diagonal selection and Fatou's lemma, a subsequence $(w_n)$ converges coordinatewise to some $x\in\ell_1(\K)$. Inductively choose strictly increasing sequences $(n_k)$ and $(m_k)$, starting with $m_0=0$, such that
\[
 \sum_{i\le m_{k-1}}|(w_{n_k}-x)_i|<2^{-k},\qquad
 \sum_{i>m_k}|(w_{n_k}-x)_i|<2^{-k}.
\]
Let $v_k=w_{n_k}$ and let $u_k$ be the restriction of $v_k-x$ to $(m_{k-1},m_k]$. Then $\|v_k-x-u_k\|_1<2^{1-k}$. Since both $(\|u_k\|_1)$ and $(\sum_i(u_k)_i)$ are bounded, passing to a further subsequence gives the two scalar limits in \eqref{eq:blocks}.

Let $P_M$ denote the projection onto the first $M$ coordinates and $Q_M=I-P_M$. Choose $M$ so that $\|Q_Mx\|_1<\eta/3$. After deleting finitely many terms we have $P_Mu_k=0$, $\|u_k\|_1>L-\eta/3$ and $\|v_k-x-u_k\|_1<\eta/3$ for every $k$. Consequently,
\begin{align*}
 \left\|\sum_{k=1}^N\lambda_k v_k\right\|_1
 &\ge\left\|Q_M\sum_{k=1}^N\lambda_k v_k\right\|_1\\
 &\ge\sum_{k=1}^N|\lambda_k|\|u_k\|_1
 -\left|\sum_{k=1}^N\lambda_k\right|\|Q_Mx\|_1\\
 &\quad-\sum_{k=1}^N|\lambda_k|\|v_k-x-u_k\|_1\\
 &\ge (L-\eta)\sum_{k=1}^N|\lambda_k|,
\end{align*}
where we have used the disjointness of the supports of the $u_k$.
\end{proof}

\begin{proof}[Proof of Theorem \ref{thm:main}]
The necessity of $1\le a\le b\le2a+1$ follows from \eqref{eq:kalenda}. Fix $1\le a\le b\le2a+1$ and set
\[
 s=\min\{b,2a\},\qquad t=b-s.
\]
Then $a\le s\le2a$ and $0\le t\le1$. Write $\Sigma(x)=\sum_{i=1}^\infty x_i$ for $x\in\ell_1(\K)$. For $x=(x_i)_{i\in\N}\in\ell_1(\K)$, put
\begin{equation}\label{eq:norm}
 \|x\|_{a,b}=\max\left\{\|x\|_1,\ a\sup_{i,j\in\N}|x_i-x_j|,\ \sup_{j\in\N}|t\Sigma(x)+sx_j|\right\}.
\end{equation}
Write $X=(\ell_1(\K),\|\cdot\|_{a,b})$. The last two terms in \eqref{eq:norm} are seminorms dominated by $a\|x\|_1$ and $b\|x\|_1$, respectively. Thus $\|\cdot\|_{a,b}$ is a norm and
\begin{equation}\label{eq:equiv}
 \|x\|_1\le\|x\|_{a,b}\le b\|x\|_1.
\end{equation}
In particular, $X$ is a Banach space.

\subsection*{The strong Schur constant: upper bound}
Fix $d\in(0,2]$, $\varepsilon>0$ and a normalized $d$-separated sequence in $X$. Apply Lemma~\ref{lem:block} to obtain a subsequence $(v_k)$ and vectors $x,u_k$ satisfying \eqref{eq:blocks}. By \eqref{eq:equiv}, $\|v_k-x-u_k\|_{a,b}\to0$. For $k\ne l$, disjointness of the supports gives
\begin{align*}
 \|u_k-u_l\|_1&=\|u_k\|_1+\|u_l\|_1,\\
 a\sup_{i,j}|(u_k-u_l)_i-(u_k-u_l)_j|&\le a(\|u_k\|_1+\|u_l\|_1),\\
 \sup_j|t\Sigma(u_k-u_l)+s(u_k-u_l)_j|
 &\le t|\Sigma(u_k)-\Sigma(u_l)|\\
 &\quad+s\max\{\|u_k\|_1,\|u_l\|_1\}.
\end{align*}
As $a\ge1$ and $s\le2a$, we obtain
\[
 d\le\ca_X((v_k)_k)\le\max\{2L,2aL,sL\}=2aL.
\]
Choose $\eta>0$ with $L-\eta>(2a/d+\varepsilon)^{-1}$. By \eqref{eq:blocklower}, after deleting finitely many terms,
\[
 \left\|\sum_{k=1}^N\lambda_k v_k\right\|_{a,b}
 \ge\left\|\sum_{k=1}^N\lambda_k v_k\right\|_1
 \ge\frac{1}{2a/d+\varepsilon}\sum_{k=1}^N|\lambda_k|
\]
for all $N\in\N$ and $\lambda_1,\ldots,\lambda_N\in\K$. The upper estimate in \eqref{eq:strong} follows from normalization. Hence $\StrSch(X)\le a$ and, by \cite[Proposition 5.6]{K}, $\SubSch(X)\le a$.

\subsection*{The Schur constant: upper bound}
Denote by $\delta_1$ and $\ca_1$ the quantities computed with the usual norm of $\ell_1(\K)$. By \cite[Example 5.9(1)]{K} and \cite[Theorem 1.3]{KS}, $\ca_1(z)=\delta_1(z)$ for every $z\in\ell_\infty(\ell_1(\K))$. From \eqref{eq:equiv}, the usual dual unit ball is contained in $B_{X^*}$, so $\delta_1(z)\le\delta(z)$. Hence
\begin{equation}\label{eq:Bupper}
 \ca(z)\le b\ca_1(z)=b\delta_1(z)\le b\delta(z),
\end{equation}
and $\Sch(X)\le b$.

\subsection*{Optimality of the constants}
Let $(e_n)_n$ be the canonical basis of $\ell_1(\K)$. Equation \eqref{eq:norm} gives, for $n\ne m$,
\begin{equation}\label{eq:basis}
 \|e_n\|_{a,b}=b,\qquad \|e_n-e_m\|_{a,b}=\max\{2,2a,s\}=2a.
\end{equation}
Let $f\in B_{X^*}$. The sequence $(f(e_n))_n$ is bounded by $b$. For any subsequence $(e_{n_k})_k$ such that $f(e_{n_k})\to\lambda\in\K$, we have
\[
 \left|\frac1N\sum_{k=1}^N f(e_{n_k})\right|
 \le\left\|\frac1N\sum_{k=1}^Ne_{n_k}\right\|_{a,b}
 =\max\{1,a/N,t+s/N\}.
\]
Letting $N\to\infty$ yields $|\lambda|\le1$. Every bounded scalar sequence has a convergent subsequence along which its moduli approach their limit superior. Hence $\limsup_n|f(e_n)|\le1$ for every $f\in B_{X^*}$. Since $\ca_\K((f(e_n))_n)\le 2\limsup_n|f(e_n)|$, we obtain $\delta((e_n)_n)\le2$. For the reverse inequality, the functional $g:X\to\K$ defined by $g(z)=\sum_{i=1}^\infty(-1)^iz_i$ belongs to $B_{X^*}$ and $\ca_\K((g(e_n))_n)=2$. For every subsequence $(e_{n_k})_{k\in\N}$ the same argument works with signs alternating along the indices $n_k$. Hence
\[
 \delta((e_n)_n)=\tdelta((e_n)_n)=2,\qquad \tca((e_n)_n)=2a.
\]
It follows that $a\le\SubSch(X)\le\StrSch(X)\le a$, so both constants are equal to $a$.

Consider the sequence
\[
 y=(e_1,-e_1,e_2,-e_2,\ldots).
\]
The bound on $\limsup_n|f(e_n)|$ implies $\delta(y)\le2$. The functional $\Sigma\in B_{X^*}$ gives the reverse inequality, so $\delta(y)=2$. By \eqref{eq:basis}, $\ca(y)=2b$, since each pair $e_n,-e_n$ has distance $2b$ and all terms have norm $b$. Therefore $\Sch(X)=b$.

\subsection*{Banach--Mazur distance}
By \eqref{eq:equiv}, the identity isomorphism gives
\[
 d(X,\ell_1(\K))\le b.
\]
Let $T:X\to\ell_1(\K)$ be any bounded linear isomorphism. For $z=(z_n)_{n\in\N}\in\ell_\infty(X)$, write $Tz=(Tz_n)_{n\in\N}$. By the $1$-Schur property of $\ell_1(\K)$,
\[
 \ca(z)\le\|T^{-1}\|\ca_1(Tz)
 =\|T^{-1}\|\delta_1(Tz)
 \le\|T^{-1}\|\|T\|\delta(z).
\]
Thus $b=\Sch(X)\le\|T\|\|T^{-1}\|$. Taking the infimum over $T$ gives $b\le d(X,\ell_1(\K))$ and completes the proof.
\end{proof}

\begin{corollary}\label{cor:sharp}
For each $c\ge1$ and either scalar field, $2c+1$ is the smallest constant $b$ for which the implication
\[
 \tca_X(z)\le c\delta_X(z)\ \text{ for all }z\in\ell_\infty(X)
 \quad\Longrightarrow\quad X\text{ is }b\text{-Schur}
\]
holds for every Banach space $X$. The same optimal constant holds when the hypothesis is replaced by the $c$-strong Schur property.
\end{corollary}
\begin{proof}
The upper bounds follow from \cite[Propositions 5.1 and 5.6]{K}. Theorem~\ref{thm:main}, applied to $(a,b)=(c,2c+1)$, gives a Banach space with $\SubSch(X)=\StrSch(X)=c$ and $\Sch(X)=2c+1$.
\end{proof}

\section{A question on the \texorpdfstring{$1$}{1}-Schur property}
\begin{proposition}\label{prop:counter}
Let $1<p<\infty$. The Banach space $X=\K\oplus_p\ell_1(\K)$ has both the $1$-Schur and $1$-strong Schur properties but fails the condition of \cite[Proposition 6.4]{K}.
\end{proposition}
\begin{proof}
The norm of $(\alpha,u)\in X$ is
\[
 \|(\alpha,u)\|_X=(|\alpha|^p+\|u\|_1^p)^{1/p}.
\]
The finite-dimensional space $\K$ is $1$-Schur, and so is $\ell_1(\K)$ by \cite[Example 5.9(1)]{K} and \cite[Theorem 1.3]{KS}. The $\ell_p$-norm on $\mathbb R^2$ lies between the $\ell_\infty$- and $\ell_1$-norms. Therefore $X$ is $1$-Schur by \cite[Theorem 7.1]{K}.

Fix $d\in(0,2]$, $\varepsilon>0$ and a normalized $d$-separated sequence $((\alpha_n,v_n))$ in $X$. Passing to a subsequence, we may assume that $(\alpha_n)$ converges and, after relabelling, that $(v_n)$ satisfies \eqref{eq:blocks}. Then, for $n\ne m$,
\[
 d\le\big(|\alpha_n-\alpha_m|^p+\|v_n-v_m\|_1^p\big)^{1/p}\longrightarrow2L
\]
as $n,m\to\infty$, by disjointness of the blocks. Choose $\eta>0$ with $L-\eta>(2/d+\varepsilon)^{-1}$. Lemma~\ref{lem:block} gives, after deleting finitely many terms,
\begin{align*}
 \left\|\sum_{n=1}^N\lambda_n(\alpha_n,v_n)\right\|_X
 &\ge\left\|\sum_{n=1}^N\lambda_n v_n\right\|_1\\
 &\ge\frac{1}{2/d+\varepsilon}\sum_{n=1}^N|\lambda_n|
\end{align*}
for all $N\in\N$ and $\lambda_1,\ldots,\lambda_N\in\K$. The upper estimate in \eqref{eq:strong} follows from normalization. Thus $X$ has the $1$-strong Schur property.

Let $(e_n)_{n\in\N}$ be the canonical basis of $\ell_1(\K)$ and set $z_n=(1,e_n)$ for every $n\in\N$. Suppose that a subsequence $(z_{n_k})_{k\in\N}$ and a vector $x=(\alpha,u)\in X$ satisfy \eqref{eq:add}. Writing $u=(u_i)_{i\in\N}$, we have $u_{n_k}\to0$ and
\[
 \|e_{n_k}-u\|_1=|1-u_{n_k}|+\|u\|_1-|u_{n_k}|\longrightarrow1+\|u\|_1.
\]
Consequently,
\begin{equation}\label{eq:L}
 L:=\lim_{k\to\infty}\|z_{n_k}-x\|_X
 =\big(|1-\alpha|^p+(1+\|u\|_1)^p\big)^{1/p}\ge1.
\end{equation}
Fix $j\in\N$. The vector $y=-z_{n_j}$ belongs to the linear span in \eqref{eq:add}. For each $l\in\N$,
\[
 \|z_{n_{j+l}}-z_{n_j}\|_X=\|e_{n_{j+l}}-e_{n_j}\|_1=2.
\]
Applying \eqref{eq:add} to this $y$ and the further subsequence $(z_{n_{j+l}})_{l\in\N}$ gives
\[
 2=\|x-z_{n_j}\|_X+L.
\]
Letting $j\to\infty$, we obtain $L=1$. By \eqref{eq:L}, this implies $\alpha=1$ and $u=0$. With $y=0$ in \eqref{eq:add}, we then have
\[
 2^{1/p}=\|x\|_X+L=2,
\]
which is impossible for $p>1$.
\end{proof}

\section*{Funding}
D. N\'u\~nez-Alarc\'on and D. Pellegrino are supported by the Conselho Nacional de Desenvolvimento Cient\'ifico e Tecnol\'ogico (CNPq, Brazil) through Grants 406457/2023-9 and 403964/2024-5 (CNPq/MCTI Calls 10/2023 and 16/2024, respectively). D. Pellegrino is also supported by CNPq Grant 305807/2025-0.

\section*{Competing interests}
The authors declare that they have no competing interests.

\section*{Data availability}
No datasets were generated or analyzed in this study.

\section*{Use of generative AI}
The authors used ChatGPT (OpenAI) to assist in the development and implementation of the solution strategy, particularly in the choice of parameters. The authors take full responsibility for the mathematical content of the article.


\end{document}